\documentclass[runningheads]{llncs}
\usepackage{amsmath}
\usepackage{amsfonts}
\usepackage{amssymb}
\usepackage{enumerate}
\usepackage{xcolor}

\newcommand\Order{\mathcal{O}}
\newcommand\nA{\mathcal{A}}
\newcommand\nS{\mathcal{S}}

\newcommand\nB{\mathcal{B}}
\newcommand\nW{\mathcal{W}}
\newcommand\nE{\mathcal{E}}
\newcommand\nL{\mathcal{L}}

\newcommand{\ee}{\mathrm{e}}

\renewcommand{\AA}{\mathtt{A}}
\newcommand{\BB}{\mathtt{B}}
\renewcommand{\gg}{\mathfrak{g}}
\newcommand{\Id}{\mathrm{Id}}
\newcommand{\grade}{\mathrm{grade}}
\newcommand{\coeff}{\mathrm{coeff}}

\begin{document}
%
%\title{Contribution Title\thanks{Supported by organization x.}}
\title{An  Algorithm for Computing Coefficients of Words in Expressions
Involving Exponentials and its Application to the
%Generation of Order Conditions for
Construction of
Exponential Integrators}
\titlerunning{Computing Coefficients of Words}
\author{Harald Hofst\"atter\inst{1}\orcidID{0000-0003-0638-9611} \and
Winfried Auzinger\inst{2}\orcidID{0000-0002-9631-2601} \and
Othmar Koch\inst{3}\orcidID{0000-0002-1956-221X}
}
\authorrunning{H. Hofst\"atter et al.}
\institute{Reitschachersiedlung 4/6, A--7100 Neusiedl am See, Austria
\email{hofi@harald-hofstaetter.at}
\url{http://harald-hofstaetter.at} \and
Vienna University of Technology, Institute of Analysis and Scientific Computing, Wiedner
Hauptstra{\ss}e 8--10, A-1040 Wien, Austria
\email{w.auzinger@tuwien.ac.at}
\url{http://asc.tuwien.ac.at/\~{}winfried} \and
University of Vienna, Institute of Mathematics, Oskar-Morgenstern-Platz 1, A-1090 Wien, Austria
\email{othmar@othmar-koch.org}
\url{http://www.othmar-koch.org}}
\maketitle              % typeset the header of the contribution
\begin{abstract}
This paper discusses an efficient implementation of the generation of order conditions
for the construction of exponential integrators  like exponential splitting and
Magnus-type methods in the computer algebra system Maple.
At the core of this implementation is a new algorithm for the computation of
coefficients of words in the formal expansion of %the expression representing
the local
error of the integrator.
The underlying theoretical background including an analysis of the structure of the
local error is briefly reviewed.
As an application the coefficients of all 8th order self-adjoint
commutator-free Magnus-type integrators involving the minimum number of 8 exponentials
are computed.
\end{abstract}

\section{Introduction}\label{Sect:Introduction}
%For the design
In the construction
of  integration schemes for the numerical solution of  evolution
equations the coefficients of the schemes are usually determined as  solutions
of certain systems of polynomial equations. The obvious requirement
that all terms up to a certain order $O(\tau^p)$
 in the Taylor
expansion of the local error with respect to the step-size $\tau$ vanish,
usually leads, if applied in a naive way, to a far over-determined system of equations.
In some cases, however, due to the special structure of the local error,
the fact that a small subset of terms vanishes already implies
that all terms up to order $O(\tau^p)$ in the Taylor
expansion vanish. This leads to a minimal, non-redundant system of equations, the so-called {\em order conditions},
see for instance \cite{auzingeretal13c}.

As it turns out, this is in particular the case for % let us consider
(generalized) exponential splitting methods
for the numerical solution of evolution equations of the form\footnote{For our
considerations, it is sufficient to discuss linear problems. The algebraic structure
underlying method construction is the same for nonlinear problems due to the calculus
of Lie derivatives~\cite[Section III.5.1]{haireretal02b}.}
\begin{equation*}
\partial_t u(t) = Au(t)+Bu(t),\quad t\geq t_0,\quad u(t_0)=u_0,\quad A, B\in\mathbb{C}^{d\times d}
\end{equation*}
(and also for Magnus-type integrators considered below).
Here one step
$u_{n+1}=\nS(\tau)u_n$
with step-size $\tau$ is specified by an approximation $\nS(\tau)=\nS_J(\tau)\cdots\nS_1(\tau)$
%\begin{equation*}%\label{eq:exactsol_split}
%\nS(\tau) = \ee^{\nF_{J}(\tau)}\cdots\ee^{\nF_{1}(\tau)}%\approx \nE(\tau)=\ee^{\tau(A+B)}
%\end{equation*}
of the exact solution operator
$\nE(\tau)=\ee^{\tau(A+B)}$, where the factors $\nS_j(\tau)$ are exponentials whose applications $\nS_j(\tau)y$ to vectors $y\in\mathbb{C}^d$ can be effectively computed. %with exponents $\nF_j$  being linear  combinations of $\tau A$, $\tau B$ and (iterated) commutators of $\tau A$ and $\tau B$. This includes the special case
%of {\em classical} splitting methods where either
%$\nF_j(\tau)=c_j\tau  A$
%or $\nF_j(\tau)=c_j\tau  B$, $c_j\in\mathbb{C}$.
Prototypical examples are for instance the classical
second order Strang splitting
\begin{equation}\label{eq:strang_tau}
\nS(\tau)=\ee^{\frac{1}{2}\tau B}\,\ee^{\tau A}\,\ee^{\frac{1}{2}\tau B},
\end{equation}
or the
4th order generalized splitting
\begin{equation}\label{eq:gen_split_tau}
\nS(\tau)=\ee^{\frac{1}{6}\tau B}\,\ee^{\frac{1}{2}\tau A}\,\ee^{\frac{2}{3}\tau B+\frac{1}{72}\tau^3[B,[A,B]]}
	\,\ee^{\frac{1}{2}\tau A}\,\ee^{\frac{1}{6}\tau B}
\end{equation}	
proposed in \cite{chin97,suzuki95}.

The analysis of  the structure of the local error
$$
\nL(\tau)= \nS(\tau)-\nE(\tau)=\nS_J(\tau)\cdots\nS_1(\tau) -\ee^{\tau(A+B)}
$$
and the resulting derivation of order conditions is advantageously
carried out in a purely formal way
by introducing non-commutative symbols $\AA$, $\BB$ representing respectively
 $\tau A$, $\tau B$, and considering the formal expression corresponding to
 $\nL(\tau)$ with
  $\tau A$, $\tau B$ substituted by $\AA$, $\BB$.
Thus, with a slight generalization antici\-pating an application to Magnus-type integrators,  we study  expressions of the  form
\begin{equation}\label{eq:Xintro}
X%=S-E
=\ee^{\Phi_J}\cdots\ee^{\Phi_1}-\ee^{\Omega},
\end{equation}
where  $\Phi_1,\dots,\Phi_J$, and $\Omega$ are linear combinations of
non-commutative symbols and commutators thereof.
 In particular,   in
Theorem~\ref{Thm:LeadingErrorTerm} of Section~\ref{Sect:OrderConditions} the structure of the leading term in the
formal series expansion of such  expressions is characterized,
which leads to the derivation of order conditions in
 Theorem~\ref{Thm:order_conditions}.
These theoretical considerations of Section~\ref{Sect:OrderConditions}
are essentially a review of   the theory developed in \cite[Section~2]{hofiopmag}, which extends
and generalizes results from \cite{auzingeretal13c}.

The main focus of  the present paper is on a concrete implementation of
the generation of order conditions according to Theorem~\ref{Thm:order_conditions}
in the computer algebra
system Maple.\footnote{We have used Maple 18, Maple is a trademark of Waterloo Maple Inc.} This can be realized in a very efficient way,
utilizing a new algorithm derived in \cite{hofiopmag} for the computation of coefficients
of words (i.e., finite products of  non-commutative symbols) in expressions like (\ref{eq:Xintro}) involving exponentials,
whose concise Maple implementation is reproduced in its entirety in Section~\ref{Sect:AlgoCoeffWords}.\footnote{All Maple code discussed in this paper
is also provided by the package {\tt Expocon} available at \cite{expocon}.
Additionally, this package includes routines for the generation of Lyndon words
and Lyndon bases. For simplicity, such words and bases have always been hardcoded
whenever needed in the code examples of this paper.
}
Our approach is thus a relevant contribution compared
to previous work on the generation of order conditions, see, e.g.,
\cite[Section~III.5]{haireretal02b} or \cite{casc2016}.

Finally in Section~\ref{Sect:Magnus},
after a brief review of material from~\cite{alvfeh11,hofiopmag} on
Magnus-type integrators for the numerical solution of non-autonomous
 evolution equations of the form
\begin{equation*}%\label{eq:non_auto_evolution_eq}
\partial_t u(t) = A(t)u(t),\quad t\geq t_0,\quad u(t_0)=u_0,\quad A(t)\in\mathbb{C}^{d\times d},
\end{equation*}
%some material from  \cite{alvfeh11} and \cite{hofiopmag} about
we consider  a  non-trivial  application of
% our theoretical approach
the theory
 of Section~\ref{Sect:OrderConditions}: We compute
in a systematic way the coefficients of all 8th order commutator-free self-adjoint Magnus-type methods involving the minimum number of 8 exponentials.

\section{Coefficients of Words in Expressions Involving Exponentials}
\label{Sect:AlgoCoeffWords}
Let $\nA$ denote a fixed set of non-commutative variables.
%Specifically,
%we consider
%\begin{equation}
%\nA=\{\AA,\BB\}\ \ \mbox{or}\ \ \nA=\{\AA_1,\AA_2,\ldots,\AA_K\}\ \mbox{for some $K%\geq 2$}.
%\end{equation}
Given an expression $X$ in these variables involving exponentials like (\ref{eq:Xintro})
%, e.g.,
%\begin{equation*}%\label{eq:strang}
%X=\ee^{\frac{1}{2}\BB}\ee^{\AA}\ee^{\frac{1}{2}\BB}- \ee^{\AA+\BB}\quad(\nA=\{\AA,\BB%\}),
%\end{equation*}
%\begin{equation*}%\label{eq:gen_split}
%X=\ee^{\frac{1}{6}\BB}\ee^{\frac{1}{2}\AA}\ee^{\frac{2}{3}\BB+\frac{1}{72}[\BB,[\AA,%\BB]]}\ee^{\frac{1}{2}\AA}\ee^{\frac{1}{6}\BB}- \ee^{\AA+\BB}\quad(\nA=\{\AA,\BB\}),
%\end{equation*}
%or
%\begin{equation*}%\label{eq:cf4}
%X=\ee^{\frac{1}{2}\AA_1+\frac{1}{3}\AA_2}\,\ee^{\frac{1}{2}\AA_1-\frac{1}{3}\AA_2}
%-\ee^{\Omega}\quad(\nA=\{\AA_1,\AA_2,\AA_3\})
%\end{equation*}
%with
%\begin{equation*}%\label{eq:Omega}
%\Omega=\AA_1-\tfrac{1}{6}[\AA_1,\AA_2]+\tfrac{1}{60}[\AA_1,[\AA_1,\AA_3]]-\tfrac{1}%{60}[\AA_2,[\AA_1,\AA_2]]
%	+\tfrac{1}{360}[\AA_1,[\AA_1,[\AA_1,\AA_2]]]-\tfrac{1}{30}[\AA_2,\AA_3],
%\end{equation*}
we want to calculate real or complex coefficients
$$
c_w=\coeff(w, X),\quad w\in\nA^*
$$
in
the formal expansion
\begin{equation*}
  X=\sum_{w\in\nA^*}c_w w\in\mathbb{C}\langle\langle\nA\rangle\rangle.
\end{equation*}
$X$ is thus represented as an element of
$\mathbb{C}\langle\langle\nA\rangle\rangle$, the algebra of formal power
series in the non-commutative variables in $\nA$.
Here, $\nA^*$ denotes the set of all words over the alphabet $\nA$, i.e., the
set of all finite products
(including the empty product $\Id$) of elements of $\nA$.

\subsection{A Family of Homomorphisms}
In \cite{hofiopmag} an efficient  algorithm for the computation of
$ \coeff(w,X)$ was derived, which
is based on a suitably constructed family of maps $\{\varphi_w:w\in\nA^*\}$, where for each word $w=w_1\cdots w_{\ell(w)}\in\nA^*$ of length $\ell(w)\geq 1$,
$\varphi_w(X)$ is an upper triangular matrix in
$\mathbb{C}^{(\ell(w)+1)\times(\ell(w)+1)}$ whose entries are coefficients
of subwords
of $w$ in $X$,
\begin{equation}\label{eq:phi_w}
%\varphi_w: \mathbb{C}\langle\langle{\nA}\rangle\rangle\to\mathbb{C}^{(\ell(w)+1)
%\times(\ell(w)+1)},\
%
\varphi_w(X)_{i,j} =
%\varphi_{i,j}^{(w)}), \
%\varphi_{i,j}^{(w)}=
\left\{\begin{array}{ll}
\coeff(w_{i:j-1},X),&\mbox{\quad if $i<j$,}\\
\coeff(\Id,X),&\mbox{\quad if $i=j$,}\\
0,&\mbox{\quad if $i>j$.}
\end{array}\right.
\end{equation}
Here $w_{i:j-1}=w_iw_{i+1}\cdots w_{j-1}$ denotes the subword of $w$ of length $j-i$, starting at position $i$ and ending at position $j-1$.
\begin{theorem}[{\cite[Theorem~2.4]{hofiopmag}}]\label{thm:III} The map $\varphi_w$
defined by {\rm (\ref{eq:phi_w})}
is an algebra homomorphism $$\mathbb{C}\langle\langle{\nA}\rangle\rangle\to\mathbb{C}^{(\ell(w)+1)\times(\ell(w)+1)},$$ i.e.,
\begin{enumerate}[(i)]
\item $\varphi_w$ is linear,
\begin{equation*}
\varphi_w(\alpha X+\beta Y)=\alpha\varphi_w(X)+\beta\varphi_w(Y),\quad
X,Y\in \mathbb{C}\langle\langle{\nA}\rangle\rangle,\ \alpha,\beta\in\mathbb{C};
\end{equation*}
\item $\varphi_w$ preserves the multiplicative structure,
\begin{equation*}
\varphi_w(X\cdot Y)=\varphi_w(X)\cdot\varphi_w(Y),\quad X,Y\in\mathbb{C}\langle\langle{\nA}\rangle\rangle.
\end{equation*}
\end{enumerate}

Furthermore,
if $\coeff(\Id, X)=0$, then
\begin{equation*}%\label{eq:exp_phi}
\varphi_w(\exp X) = \exp\varphi_w(X),
\end{equation*}
where the exponential of the strictly upper triangular and thus nilpotent matrix $\varphi_w(X)$ is exactly computable in
a finite number of steps.
\end{theorem}
\subsection{Maple Implementation of the Algorithm}
It follows that
for a given expression $X$, a recursive application of $\varphi_w$ (the recursion terminates with well-defined values $\varphi_w(a)$ for the ``atoms'' $a\in\nA$) yields $\varphi_w(X)$, from which one can read off $\coeff(w,X)$ as the element at the upper right corner,
$$\coeff(w,X)=\varphi_w(X)_{1,\ell(w)+1},$$
cf.~(\ref{eq:phi_w}).
By organizing this calculation in a more efficient way,
the  function \verb|phiv| defined in the Maple code displayed below computes
\begin{equation*}
  \verb|phiv|(w, X, v) = \varphi_w(X)\cdot v
\end{equation*}
for a vector $v\in\mathbb{C}^{\ell(w)+1}$
without explicitly generating the matrix $\varphi_w(X)$.
%The function \verb|phiv|
It recursively traverses the expression tree representing the expression $X$. At each node of the tree the evaluation branches out
depending on whether the current node represents
\begin{itemize}
\item
a non-commutative symbol (the atomic case which terminates the recursion),
\item
a sum of subexpressions,
\item
a product of subexpressions,
\item
a power of a subexpression,
\item
a commutator of subexpressions, or
\item
an exponential of a subexpression.
\end{itemize}
Finally,  the function \verb|wcoeff|\footnote{The function was called {\tt wcoeff} because {\tt coeff} is already defined in Maple.} computes $\coeff(w, X)$ via
\begin{equation*}
\coeff(w,X)=\mbox{first component of}\ %\varphi_w(X)\cdot(0,\dots,0,1)^T
\verb|phiv|(w, X,(0,\dots,0,1)^T ).
\end{equation*}
The elements of the alphabet $\nA$ are
represented within Maple as  non-commutative symbols, which are provided
by the package  \verb|Physics|. Note that
except for providing such non-commutative symbols (and the type \verb|Commutator|)
we do not need or use any further feature of the package \verb|Physics|.
Words $w\in\nA^*$ are represented as lists of non-commutative symbols.
\begin{verbatim}
> with(Physics):
> phiv := proc (w, X, v)
    local i, v1, v2, f, zero;
    if type(X, name) and type(X, noncommutative) then
        return [seq(`if`(op(i,w)=X, v[i+1], 0), i=1..nops(w)), 0]
    elif type(X, `+`) then
        return add(phiv(w, op(i, X), v), i=1..nops(X))
    elif type(X, `*`) then
        v1 := v; zero := [0$nops(w)+1];
        for i from nops(X) to 1 by -1 do
            v1 := phiv(w, op(i, X), v1);
            if v1=zero then return zero end if;
        end do;
        return v1
    elif type(X, anything^integer) then
        v1 := v; zero := [0$nops(w)+1];
        for i from 1 to op(2, X) do
            v1 := phiv(w, op(1, X), v1);
            if v1=zero then return zero end if;
        end do;
        return v1
    elif (type(X, function) and
          op(0, X) = Physics[Commutator]) then
        return phiv(w, op(1, X), phiv(w, op(2, X), v))
             - phiv(w, op(2, X), phiv(w, op(1, X), v))
    elif type(X, exp(anything)) then
        v1 := v; v2 := v; zero := [0$nops(w)+1]; f := 1;
        for i from 1 to nops(w) do
            f := f*i; v1 := phiv(w, op(X), v1);
            if v1=zero then return v2 end if;
            v2 := v2 + v1/f;
        end do;
        return v2
    end if;
    return [seq(X*x, x=v)]
  end proc:
> wcoeff := proc (w, X)
    return phiv(w, X, [0$nops(w), 1])[1]
  end proc:
\end{verbatim}

\section{Order Conditions for Exponential Integrators}
\label{Sect:OrderConditions}

In this section we review the theory developed in \cite[Section~2]{hofiopmag}, which extends
and generalizes results from \cite{auzingeretal13c}.
We consider expressions of the form (\ref{eq:Xintro}),
\begin{equation}\label{eq:X_gen}
X=\ee^{\Phi_J}\cdots\ee^{\Phi_1}-\ee^{\Omega},
\end{equation}
where the exponents $\Phi_1,\dots,\Phi_J$, and $\Omega$ are
linear combinations of non-commutative symbols and commutators thereof, i.e.,
 elements of
$[\mathbb{C}\langle\nA\rangle]$, the free Lie algebra generated by the non-commutative symbols of a given alphabet $\nA$, which in a natural way is embedded in
the algebra $\mathbb{C}\langle\langle\nA\rangle\rangle$  of formal power
series in these symbols.

In the applications we are interested in,
$S=\ee^{\Phi_J}\cdots\ee^{\Phi_1}$ represents an exponential integrator
for the numerical solution of an evolution equation, and $E=\ee^{\Omega}$ represents
the exact local solution operator for this equation.
We can interpret (\ref{eq:X_gen}) as the error of the approximation
$S$ of $E$, i.e.,
$X$ represents the {\em local error} of the exponential integrator $S$.

\subsection{Grading of Words and Homogeneous Lie Elements}
\label{SubSec:Grading}
We consider a grading function on the alphabet $\nA$,
\begin{equation}\label{eq:grading_symbols}
\grade(a) \in\{1,2,\dots\}%\ \mbox{given for all}
, \quad a\in\nA,
\end{equation}
and extend it to words $w=w_1\dots w_{\ell(w)}\in\nA^*$ by
$$\grade(w) = \sum_{j=1}^{\ell(w)}\grade(w_j).$$
We call $\Psi\in[\mathbb{C}\langle\nA\rangle]$ a homogeneous Lie element of grade
$q$ if it can be expanded  in
$\mathbb{C}\langle\langle\nA\rangle\rangle$
to a
linear combination of words all of the same grade $q$.
 The decomposition
\begin{equation}\label{eq:ggq}
[\mathbb{C}\langle\nA\rangle]=\bigoplus_{q=1}^{\infty}\gg_q,\quad
\gg_q=\{\mbox{homogeneous Lie elements of grade $q$}\}
\end{equation}
into a direct sum of subspaces makes $[\mathbb{C}\langle\nA\rangle]$
a graded Lie algebra, cf.~\cite{MuntheKaas957}.
\begin{remark}\label{Rmk:grading}
In the applications we are interested in,
the symbols $a\in\nA$ represent objects which depend on a (small) parameter
$\tau>0$ (for instance, a time increment). The grading (\ref{eq:grading_symbols}) is chosen such that
it reflects the order of magnitude of the represented objects,
$$
a\simeq O(\tau^{\grade(a)}),\quad a\in\nA.
$$
For example, in the case of an application to splitting methods with step-size $\tau$,
 $$\nA=\{\AA,\BB\},\quad\AA\simeq \tau A=O(\tau),\ \BB\simeq \tau B=O(\tau)
\ \Rightarrow \ \grade(\AA)=\grade(\BB)=1,
 $$
cf.~Section~\ref{Sect:Introduction}.

\end{remark}

\subsection{Leading Error Term}% of an Approximation}
 The following
theorem states that the leading error term $\Theta$ of the
 approximation  $\ee^{\Phi_J}\cdots\ee^{\Phi_1}$ of $\ee^{\Omega}$
is a homogeneous Lie element of some grade $q$.
%This grade $q$ determines the order of the approximation.
%
\begin{theorem}[{\cite[Theorem~2.1]{hofiopmag}}]\label{Thm:LeadingErrorTerm}
If for $\Phi_1,\dots,\Phi_J,\Omega\in[\mathbb{C}\langle\nA\rangle]$
the expression $\ee^{\Phi_J}\cdots\ee^{\Phi_1}-\ee^{\Omega}$ is expanded
in $\mathbb{C}\langle\langle\nA\rangle\rangle$ as
\begin{equation*}
\ee^{\Phi_J}\cdots\ee^{\Phi_1}-\ee^{\Omega}
=\sum_{w\in\nA^*}c_w w = \Theta +R ,
\end{equation*}
where
\begin{equation}\label{eq:Theta}
\Theta = \sum_{\grade(w)=q_{\mathrm{min}}}\!\!c_w w,\qquad
q_{\mathrm{min}}=\mathrm{min}\{\grade(w):w\in\nA^*, c_w\neq 0\}
\end{equation}
(and the remainder $R$ contains the terms of grade $>q_{\mathrm{min}}$),
then $\Theta$ can be represented as a linear combination of commutators, i.e.,
$\Theta$ is a homogeneous Lie element of grade $q_\mathrm{min}$.
\end{theorem}
To illustrate Theorem~\ref{Thm:LeadingErrorTerm} we consider as an example
$$X=\ee^{\frac{1}{2}\BB}\ee^{\AA}\ee^{\frac{1}{2}\BB}-\ee^{\AA+\BB},\quad
\nA=\{\AA,\BB\},\quad  \grade(\AA)=\grade(\BB)=1,$$
cf.~(\ref{eq:strang_tau}).
The following Maple code computes the coefficients
of all words of length $\leq 3$ (i.e., of all words $w$ with $\grade(w)\leq 3)$
in $X$.
\begin{verbatim}
> Physics[Setup](noncommutativeprefix = {A, B}):
> X := exp((1/2)*B)*exp(A)*exp((1/2)*B)-exp(A+B):
> W := [[A], [B], [A, A], [A, B], [B, A], [B, B],
       [A, A, A], [A, A, B], [A, B, A], [A, B, B],
       [B, A, A], [B, A, B], [B, B, A], [B, B, B]]:
> seq(wcoeff(w, X), w in W);
\end{verbatim}
%                             1   -1  -1  1   1   -1
%        0, 0, 0, 0, 0, 0, 0, --, --, --, --, --, --, 0
%                             12  6   24  12  12  24
$$0,0,0,0,0,0,0,\frac{1}{12},\frac{-1}{6},\frac{-1}{24},\frac{1}{12},
\frac{1}{12},\frac{-1}{24},0$$
It follows
\begin{align*}
\ee^{\frac{1}{2}\BB}\ee^{\AA}\ee^{\frac{1}{2}\BB}-\ee^{\AA+\BB}
&=\tfrac{1}{12}\AA\AA\BB
-\tfrac{1}{6} \AA\BB\AA
-\tfrac{1}{24}\AA\BB\BB
+\tfrac{1}{12}\BB\AA\AA
+\tfrac{1}{12}\BB\AA\BB
-\tfrac{1}{24}\BB\BB\AA+\dots\\
&=\tfrac{1}{12}[\AA,[\AA,\BB]]-\tfrac{1}{24}[[\AA,\BB],\BB]+\dots.
\end{align*}
Here the leading error term is indeed a homogeneous Lie element of grade 3.

\subsection{Symmetry}\label{subsec:symmetry}
A product of exponentials $S=\ee^{\Phi_J}\cdots\ee^{\Phi_1}$,
$\Phi_j\in[\mathbb{C}\langle\nA\rangle]$ is called
self-adjoint or symmetric,\footnote{For an exponential integrator
represented by $S$
this definition conforms with the usual
definition of a self-adjoint integrator,
e.g., $\nS(-\tau)\nS(\tau)=\Id$ for a generalized splitting method
where $\nS(\tau)$ is $S$ with $\tau A,\tau B$ substituted for $\AA,\BB$.} if
$$\Phi_{J-j+1}=\sum_k(-1)^{k+1} X_{j,k},\quad j=1,\dots,J$$
holds, where the $X_{j,k}\in\gg_{k}$ are the components
of $\Phi_j=\sum_k X_{j,k}$
%in
with respect to
the decomposition (\ref{eq:ggq}).
It follows that a single exponential $\ee^\Phi$ is self-adjoint, if and only if
$\Phi$ is a sum of homogeneous Lie elements of odd grade,
$\Phi = X_1+X_3+ \dots,\quad X_k\in\gg_k$.
It was proved in
\cite[Theorem~2.2]{hofiopmag} that  in (\ref{eq:Theta}) the grade $q_{\mathrm{min}}$
of the homogeneous Lie element $\Theta$ is necessarily odd if  $\ee^{\Phi_J}\cdots\ee^{\Phi_1}$
and $\ee^{\Omega}$ are both self-adjoint.

\subsection{Lyndon words and Lyndon bases}
For a homogeneous Lie element $\Theta$ of grade $q$  like the one given in
(\ref{eq:Theta}) let
\begin{equation}\label{eq:ThetaB}
\Theta=\sum_{b\in\nB_q}c_b\,b
\end{equation}
be its representation
in a basis $\nB_q$  of the subspace $\gg_q$ of (\ref{eq:ggq}).
Furthermore, let
$\nW_q\subset\nA^*$ be a set of words of grade $q$ such that the
matrix
\begin{equation}\label{eq:Tq}
T_q=(\coeff(w,q))_{w\in\nW_q,b\in\nB_q}
\end{equation}
is invertible. Then it follows from
$c_w=\coeff(w,\Theta)=\sum_{b\in\nB_q}c_b\,\coeff(w,b)$ for $w\in\nW_q$,
and thus $(c_w)_{w\in\nW_q}=T_q\cdot(c_b)_{b\in\nB_q}  $,
 that
the coefficients $c_b$ in (\ref{eq:ThetaB}) can be computed as
\begin{equation}\label{eq:cb}
(c_b)_{b\in\nB_q}=T_q^{-1}\cdot(c_w)_{w\in\nW_q}.
\end{equation}
Suitable choices for  such a set $\nW_q$ and basis $\nB_q$ are respectively the set of Lyndon words of grade $q$ and the corresponding Lyndon basis
\cite{duval88,lothaire97}, see Tables~\ref{tab:lyndon_a_b} and \ref{tab:lyndon_ak}.

\begin{table}[t]
\begin{center}
\caption{Lyndon words $\nW_q$ of grade $q$ and Lyndon basis $\nB_q$ of $\gg_q$ for
$\nA=\{\AA,\BB\}$ and $\grade(\AA)=\grade(\BB)=1$.
\label{tab:lyndon_a_b}
}
\begin{tabular}{|c|l|l|}
\hline
 $q$ & Lyndon words & Lyndon basis %& rightnormed basis
\\
\hline
\hline
1   & $\AA$, $\BB$
    & $\AA$, $\BB$
%    & $\AA$, $\BB$
\\
\hline
2	& $\AA\BB$
    & $[\AA,\BB]$
%    & $[\BB,\AA]$
\\
\hline
3   & $\AA\AA\BB$, $\AA\BB\BB$
	& $[\AA,[\AA,\BB]]$,\ $[[\AA,\BB],\BB]$
%	&  $[\AA,[\BB,\AA]]$,  $[\BB,[\BB,\AA]]$
\\
\hline
4   &  $\AA\AA\AA\BB$,  $\AA\AA\BB\BB$, $\AA\BB\BB\BB$
    &  $[\AA,[\AA,[\AA,\BB]]]$,\ $[\AA,[[\AA,\BB],\BB]]$,\ $[[[\AA,\BB],\BB],\BB]$
%    & $[\AA,[\AA,[\BB,\AA]]]$, $[\BB,[\AA,[\BB,\AA]]]$, $[\BB,[\BB,[\BB,\AA]]]$
\\
\hline
5   &$\AA\AA\AA\AA\BB$, $\AA\AA\AA\BB\BB$, $\AA\AA\BB\AA\BB$,
	&$[\AA,[\AA,[\AA,[\AA,\BB]]]]$,\ $[\AA,[\AA,[[\AA,\BB],\BB]]]$,\ $[[\AA,[\AA,\BB]],[\AA,\BB]]$,
%    & $[\AA,[\AA,[\AA,[\BB,\AA]]]]$, $[\BB,[\AA,[\AA,[\BB,\AA]]]]$,$[\AA,[\BB,[\AA,[\BB,\AA]]]]$,
    \\
    &$\AA\AA\BB\BB\BB$, $\AA\BB\AA\BB\BB$, $\AA\BB\BB\BB\BB$
    &$[\AA,[[[\AA,\BB],\BB],\BB]]$,\ $[[\AA,\BB],[[\AA,\BB],\BB]]$,\ $[[[[\AA,\BB],\BB],\BB],\BB]$
%    &$[\BB,[\BB,[\AA,[\BB,\AA]]]]$,$[\AA,[\BB,[\BB,[\BB,\AA]]]]$, $[\BB,[\BB,[\BB,[\BB,\AA]]]]$
    \\
\hline
\end{tabular}
\end{center}
\end{table}

\begin{table}[t]
\begin{center}
\caption{Lyndon words $\nW_q$ of grade $q$ and Lyndon basis $\nB_q$ of $\gg_q$ for  $\nA=\{\AA_1,\dots,\AA_q\}$ and $\grade(\AA_k)=k$.
\label{tab:lyndon_ak}
}
\begin{tabular}{|c|l|l|}
\hline
$q$ &Lyndon words & Lyndon basis %& rightnormed basis\!\!\!
\\
\hline
\hline
1 & $\AA_1$ & $\AA_1$ %& $\AA_1$
\\
\hline
2 &  $\AA_2$ & $\AA_2$ %& $\AA_2$
\\
\hline
3 & $\AA_1\AA_2$, $\AA_3$
  & $[\AA_1,\AA_2]$, $\AA_3$
%  & $[\AA_2,\AA_1]$, $\AA_3$
\\
\hline
4 & $\AA_1\AA_1\AA_2$, $\AA_1\AA_3$, $\AA_4$
  & $[\AA_1,[\AA_1,\AA_2]]$, \ $[\AA_1,\AA_3]$, \ $\AA_4$
%  &$[\AA_1,[\AA_2,\AA_1]]$, $[\AA_3,\AA_1]$, $\AA_4$
\\
\hline
5 & $\AA_1\AA_1\AA_1\AA_2$, $\AA_1\AA_1\AA_3$, $\AA_1\AA_2\AA_2$,
  & $[\AA_1,[\AA_1,[\AA_1,\AA_2]]]$, \ $[\AA_1,[\AA_1,\AA_3]]$, \ $[[\AA_1,\AA_2],\AA_2]$,
%  & $[\AA_1,[\AA_1,[\AA_2,\AA_1]]]$, $[\AA_1,[\AA_3,\AA_1]]$, $[\AA_2,[\AA_2,\AA_1]]$,
 \\
 & $\AA_1\AA_4$, $\AA_2\AA_3$, $\AA_5$
 & $[\AA_1,\AA_4]$, \ $[\AA_2,\AA_3]$, \ $\AA_5$
%  $[\AA_4,\AA_1]$,  &$[\AA_3,\AA_2]$, $\AA_5$
  \\
\hline
\end{tabular}
\end{center}
\end{table}

\subsection{Order Conditions}
The following  main result of this section is
an easy consequence of the previous considerations.
\begin{theorem}[{\cite[Theorem~2.3]{hofiopmag}}]\label{Thm:order_conditions}
If for $\Phi_1,\dots,\Phi_J,\Omega\in[\mathbb{C}\langle\nA\rangle]$ the order conditions
\begin{equation}\label{eq:order_conditions}
c_w=\coeff(w,\ee^{\Phi_J}\cdots\ee^{\Phi_1}-\ee^{\Omega})=0,\quad
w\in\bigcup_{q=1}^p\nW_q
\end{equation}
are satisfied
for all Lyndon words
 of grade $q\leq p$,
then
\begin{equation}\label{eq:order_p}
c_w=\coeff(w,\ee^{\Phi_J}\cdots\ee^{\Phi_1}-\ee^{\Omega})=0,\quad
w\in\nA^*, \ \grade(w)\leq p,
\end{equation}
and thus $q_\mathrm{min}\geq p+1$ in {\rm (\ref{eq:Theta})}.

If  $\ee^{\Phi_J}\cdots\ee^{\Phi_1}$
and $\ee^{\Omega}$ are both self-adjoint, then we may assume that $p$ is even, and {\rm (\ref{eq:order_p})} holds already
if the order conditions {\rm (\ref{eq:order_conditions})} are satisfied only for all
Lyndon words of {\em odd} grade $q\leq p$.
\end{theorem}
In view of Remark~\ref{Rmk:grading} in
Subsection~\ref{SubSec:Grading}
we can interpret (\ref{eq:order_p}) as the statement
$$\ee^{\Phi_J}\cdots\ee^{\Phi_1}-\ee^{\Omega}\simeq O(\tau^{p+1}),$$
i.e.,
$\ee^{\Phi_J}\cdots\ee^{\Phi_1}$ is an approximation of $\ee^{\Omega}$ of order $p+1$.
\subsection{Example}
For $\nA=\{\AA,\BB\}$ with $\grade(\AA)=\grade(\BB)=1$ we want to determine the
parameters $a,b,c,d\in\mathbb{R}$ such that
$S=\ee^{b\BB}\ee^{a\AA}\ee^{c\BB+d[\BB,[\AA,\BB]]}\ee^{a\AA}\ee^{b\BB}$
is a 5th order approximation of  $E=\ee^{\AA+\BB}$. Since the ansatz
$S$ and the expression $E$ are both self-adjoint in the sense of
Subsection~\ref{subsec:symmetry},
 we only have to consider   Lyndon words of  odd grade $\leq4$ (i.e., of odd length $\leq 4$ for $\grade(\AA)=\grade(\BB)=1$),
$$
\nW=\nW_1\cup\nW_3=\{\AA, \BB,  \AA\AA\BB, \AA\BB\BB\},
$$
cf.~Table~\ref{tab:lyndon_a_b}.
The order conditions (\ref{eq:order_conditions}) lead to 4 polynomial equations
in  4 variables $a,b,c,d$, for which the following Maple code computes a unique solution corresponding
to $$S=\ee^{\frac{1}{6}\BB}\ee^{\frac{1}{2}\AA}\ee^{\frac{2}{3}\BB+\frac{1}{72}[\BB,[\AA,\BB]]}\ee^{\frac{1}{2}\AA}\ee^{\frac{1}{6}\BB},$$
cf.~(\ref{eq:gen_split_tau}).
\begin{verbatim}
> Physics[Setup](noncommutativeprefix = {A, B}):
> C := Physics[Commutator]:
> X := exp(b*B)*exp(a*A)*exp(c*B+d*C(B, C(A, B)))*
       exp(a*A)*exp(b*B) - exp(A+B):
> W := [[A], [B], [A, A, B], [A, B, B]]:
> eqs := [seq(simplify(wcoeff(w, X)), w in W)];
\end{verbatim}
%
%        [                          1      2     1  2
%        [-1 + 2 a, -1 + 2 b + c, - - + 2 a  b + - a  c,
%        [                          6            2
%            1   1    2              2    ]
%          - - + - a c  + a c b + a b  - d]
%            6   2                        ]
$$eqs:=\left[-1+2a, -1+2b+c, -\frac{1}{6}+2a^2b+\frac{1}{2}a^2c,
-\frac{1}{6}+\frac{1}{2}ac^2+acb+ab^2-d
\right]$$
\begin{verbatim}
> sol := solve(eqs);
\end{verbatim}
%                 /    1      1      2      1 \
%                { a = -, b = -, c = -, d = -- }
%                 \    2      6      3      72/
$$sol:=\left\{a=\frac{1}{2},b=\frac{1}{6},c=\frac{2}{3},d=\frac{1}{72}\right\}$$
%
%\medskip
%
Next we compute the leading error term $\Theta$
of the approximation $S$ of $E=\ee^{A+B}$, cf.~(\ref{eq:Theta}).
To this end we take
$\nW_5$ and $\nB_5$ from Table~\ref{tab:lyndon_a_b} and compute  $T_5$,
$(c_w)_{w\in\nW_5}$, and $(c_b)_{b\in\nB_5}$ according
to (\ref{eq:Tq}), (\ref{eq:cb}):
\begin{verbatim}
> W5 := [[A, A, A, A, B], [A, A, A, B, B], [A, A, B, A, B],
         [A, A, B, B, B], [A, B, A, B, B], [A, B, B, B, B]]:
> B5 := [C(A, C(A, C(A, C(A, B)))), C(A, C(A, C(C(A, B), B))),
         C(C(A, C(A, B)), C(A, B)), C(A, C(C(C(A, B), B), B)),
         C(C(A, B), C(C(A, B), B)), C(C(C(C(A, B), B), B), B)]:
> T5 := Matrix([seq([seq(wcoeff(w, b), b in B5)], w in W5)]);
\end{verbatim}
%                [1     0    0     0    0    0]
%                [                            ]
%                [0     1    0     0    0    0]
%                [                            ]
%                [0    -2    1     0    0    0]
%                [                            ]
%                [0     0    0     1    0    0]
%                [                            ]
%                [0     0    0    -3    1    0]
%                [                            ]
%                [0     0    0     0    0    1]
$$
T5:=\left[
\begin{array}{rrrrrr}
1 & 0 & 0 & 0 & 0 & 0 \\
0& 1 & 0 & 0 & 0 & 0\\
0& -2 & \phantom{-}1 & 0 & \phantom{-}0 & \phantom{-}0 \\
0& 0& 0& 1 & 0 & 0 \\
0& 0& 0& -3 & 1 & 0\\
0& 0& 0& 0& 0& 1
\end{array}
\right]
$$
\begin{verbatim}
> c_w := [seq(wcoeff(w, subs(sol, X)), w in W5)];
\end{verbatim}
%             [ 1     -7    1     7    -1    -41  ]
%             [----, ----, ---, -----, ---, ------]
%             [2880  8640  480  12960  720  155520]
$$c\_w := \left[\frac{1}{2880},\frac{-7}{8640},\frac{1}{480},
\frac{7}{12960},\frac{-1}{720},\frac{-41}{155520}
\right]$$
\begin{verbatim}
> c_b :=  evalm(LinearAlgebra[MatrixInverse](T5) &* c_w);
\end{verbatim}
%            [ 1     -7    1      7     1     -41  ]
%            [----, ----, ----, -----, ----, ------]
%            [2880  8640  2160  12960  4320  155520]
$$c\_b:=\left[\frac{1}{2880},\frac{-7}{8640},\frac{1}{2160},
\frac{7}{12960},\frac{1}{4320},\frac{-41}{155520}
\right]$$
Altogether, we obtain a representation of the leading error term,
\begin{align*}
&\ee^{\frac{1}{6}\BB}\ee^{\frac{1}{2}\AA}\ee^{\frac{2}{3}\BB+\frac{1}{72}[\BB,[\AA,\BB]]}\ee^{\frac{1}{2}\AA}\ee^{\frac{1}{6}\BB}-\ee^{\AA+\BB}= \Theta + \dots\\
&\ =\tfrac{1}{2880}[\AA,[\AA,[\AA,[\AA,\BB]]]]-\tfrac{7}{8640}[\AA,[\AA,[[\AA,\BB],\BB]]]+\tfrac{1}{2160}[[\AA,[\AA,\BB]],[\AA,\BB]]\\
&\quad\ +\tfrac{7}{12960}[\AA,[[[\AA,\BB],\BB],\BB]]+\tfrac{1}{4320}[[\AA,\BB],[[\AA,\BB],\BB]]-\tfrac{41}{155520}[[[[\AA,\BB],\BB],\BB],\BB]+\dots.
\end{align*}	
Here,  the dots represent terms of grade higher than five.

\section{Magnus-Type Integrators}
\label{Sect:Magnus}
In this section we apply the theory of Section~\ref{Sect:OrderConditions}
with the aim of constructing
 Magnus-type integrators for
the numerical solution of
 non-autonomous evolution equations
\begin{equation*}%\label{eq:non_auto_evolution_eq}
\partial_t u(t) = A(t)u(t),\quad t\geq t_0,\quad u(t_0)=u_0,\quad A(t)\in\mathbb{C}^{d\times d}.
\end{equation*}
One step $(t_n,u_n)\mapsto(t_{n+1},u_{n+1})$ of step-size $\tau$ of
such an integrator is given by
\begin{equation}\label{eq:non_auto_step}
t_{n+1}=t_n+\tau, \quad u_{n+1}=\nS(\tau, t_n)u_n,
\end{equation}
where $\nS(\tau, t_n)\approx\nE(\tau, t_n)$ approximates the exact local solution operator
\begin{equation}\label{eq:exactsol}
\nE(\tau, t_n)=\ee^{\Omega(\tau,t_n)}
\end{equation}
given by the {\em Magnus} series $\Omega=\Omega(\tau, t_n)$, see
\cite[Section~IV.7]{haireretal02b}.
\subsection{Legendre expansions}
To construct $\nS(\tau,t_n)$ we
expand $A(t_n+t)$ on the interval $[t_n, t_n+\tau]$
into a series of shifted
Legendre polynomials,
\begin{equation*}%\label{eq:A_expansion}
A(t_n+t)=A_1\tilde{P}_0(t)+A_2\tilde{P}_1(t)+A_3\tilde{P}_2(t)+\dots,\quad t\in[0,\tau]
\end{equation*}
with
\begin{equation*}%\label{eq:legendre_zeugs}
\tilde P_k(t) = \frac{1}{\tau}P_{k}\left(\frac{t}{\tau}\right),\quad
P_{k}(x)=(-1)^k\sum_{j=0}^k{k \choose j}{k+j \choose j}(-1)^j x^j,
\end{equation*}
see~\cite[Section~3.1]{alvfeh11}.
The matrix-valued coefficients $A_k$ given by
\begin{equation}\label{eq:A_integral}
A_k%=(2k-1)\tau\int_0^\tau A(t_n+t)\tilde{P}_{k-1}(t)\,\mathrm{d}t
=(2k-1)\tau\int_0^1 P_{k-1}(x)A(t_n+\tau x)\,\mathrm{d}x
\end{equation}
depend on both $t_n$ and $\tau$ and satisfy
\begin{equation} \label{eq:legendre_coeffs}
A_k=\Order(\tau^k).
\end{equation}
In terms of these coefficients
the Magnus series in (\ref{eq:exactsol}) is given by
\begin{eqnarray}
	\Omega &=&  A_1-\tfrac{1}{6}[A_1,A_2]+\tfrac{1}{60}[A_1,[A_1,A_3]]-\tfrac{1}{60}[A_2,[A_1,A_2]]\nonumber\\
	&&+\tfrac{1}{360}[A_1,[A_1,[A_1,A_2]]]-\tfrac{1}{30}[A_2,A_3]%-\frac{1}{70}[A3,A4] A+
	+\dots.
	\label{eq:magnus_series}
\end{eqnarray}	
see~\cite[Section~3.2]{alvfeh11}.
\subsection{Order Conditions for Magnus-Type Integrators}
We consider Magnus-type integrators (\ref{eq:non_auto_step}) of the form
\begin{equation}\label{eq:magnus_type_int}
\nS(\tau,t_n)=\ee^{\widetilde{\Phi}_J(\tau,t_n)}\cdots\ee^{\widetilde{\Phi}_1(\tau,t_n)},
\end{equation}
where the $\widetilde{\Phi}_j$ are
linear combinations of (commutators of)
approximations
$\widetilde{A}_k\approx A_k$
obtained by applying a suitable quadrature formula to (\ref{eq:A_integral}).
To  apply the theory of
Section~\ref{Sect:OrderConditions} to such integrators, we note that
(\ref{eq:magnus_type_int}) formally corresponds to an expression
$$
S=\ee^{\Phi_J}\cdots\ee^{\Phi_1}
$$
with Lie elements $\Phi_j\in[\mathbb{C}\langle\nA\rangle]$ over
an alphabet $\nA=\{\AA_1,\AA_2,\dots\}$ with symbols $\AA_k$ representing
$\widetilde{A}_k\approx A_k$, and with
a grading
$$
\grade(\AA_k)=k
$$
corresponding to (\ref{eq:legendre_coeffs}), cf.~Remark~\ref{Rmk:grading} in
Subsection~\ref{SubSec:Grading}.

To set up  order conditions according to
Theorem~\ref{Thm:order_conditions}
we need the Lyndon words from Table~\ref{tab:lyndon_ak}, and,  furthermore,
we have to consider
$\ee^\Omega$ with $\Omega$ from (\ref{eq:magnus_series}),
%(with $A_k$ substituted by $\AA_k$),
which is self-adjoint in the sense of
Subsection~\ref{subsec:symmetry}.
For coefficients of words $w\in\nA^*$ in $\ee^\Omega$ (which
in principle could be calculated with the algorithm
of Section~\ref{Sect:AlgoCoeffWords}) we use the explicit formula
\begin{equation}\label{eq:coeff_rhs_magnus_type}
\mathrm{coeff}(\AA_{d_1}\cdots\AA_{d_\ell},\ee^{\Omega}) %=\mathrm{coeff}(w_1\dots w_{n_w},\ee^{\Omega})
=\sum_{(k_1,\dots,k_{\ell})\atop 1\leq k_l\leq d_l}
\prod_{j=1}^{\ell}\frac{
(-1)^{d_j+k_j}{d_j-1 \choose k_j-1}{d_j+k_j-2 \choose k_j-1}
}{\sum_{i=j}^{\ell}k_i}
\end{equation}
proven in \cite[Theorem~4.1]{hofiopmag}.

\subsection{8th Order Commutator-Free Magnus-Type Integrators}
In this section we construct 8th order  self-adjoint commutator-free
integrators
involving a minimum number of exponentials.
In an ansatz for such a scheme only the generators
$\AA_1,\AA_2,\AA_3,\AA_4$ have to be considered, because it can be shown that $\coeff(w, \ee^\Omega)=0$ for all words $w$ of $\grade(w)\leq 8$  containing $\AA_k$ with $k\geq 5$, see \cite[Section~3.3]{alvfeh11}.\footnote{Of course, this follows also
from (\ref{eq:coeff_rhs_magnus_type}) by a direct computation.}
 % and \cite[Section~4.1]{hofiopmag}.
Corresponding to 22 Lyndon words of odd grade $\leq 8$ over the alphabet
$\nA=\{\AA_1,\AA_2,\AA_3,\AA_4\}$ there are
22 order-conditions to be considered, see Theorem~\ref{Thm:order_conditions}. This implies an ansatz
involving 11 exponentials and 22 parameters to be determined;
 an 8th order scheme of this form  was  derived
 in \cite[Section~4.4]{alvfeh11}.
Also with this approach we found in \cite[Section~4.4]{hofiopmag}
a scheme where some exponentials commute, which can thus be joined
together, resulting in an 8th order scheme involving only 8 exponentials.
These schemes were found by a rather brute force computation. In contrast,
using the following
 Maple code we are able to
compute the coefficients of all 8th order self-adjoint schemes with 8 exponentials in a more systematic and efficient\footnote{Compared with the effort indicated in
\cite[Section~4.4]{alvfeh11}.} way.

First, we define the self-adjoint ansatz for a scheme involving 8 exponentials.
\begin{verbatim}
> Physics[Setup](noncommutativeprefix = {A}):
> S := exp(f11*A1-f12*A2+f13*A3-f14*A4)*
       exp(f21*A1-f22*A2+f23*A3-f24*A4)*
       exp(f31*A1-f32*A2+f33*A3-f34*A4)*
       exp(f41*A1-f42*A2+f43*A3-f44*A4)*
       exp(f41*A1+f42*A2+f43*A3+f44*A4)*
       exp(f31*A1+f32*A2+f33*A3+f34*A4)*
       exp(f21*A1+f22*A2+f23*A3+f24*A4)*
       exp(f11*A1+f12*A2+f13*A3+f14*A4):
\end{verbatim}
Next we set up the 8 equations corresponding to the 8 Lyndon words of odd grade $\leq 8$ involving only the generators $\AA_1, \AA_2$. The right-hand sides
of these equations were computed using (\ref{eq:coeff_rhs_magnus_type})
and are hardcoded here for simplicity.
\begin{verbatim}
> W12 := [[A1], [A1, A2], [A1, A1, A1, A2], [A1, A2, A2],
          [A1, A1, A1, A1, A1, A2], [A1, A1, A1, A2, A2],
          [A1, A1, A2, A1, A2], [A1, A2, A2, A2]]:
> rhs12 := [1, -1/6, -1/40, 1/60, -1/1008, 1/420, 1/2520, -1/840]:
> vars12 := [f11, f21, f31, f41, f12, f22, f32, f42]:
> eqs12 := [seq(expand(wcoeff(w, S)), w in W12)] - rhs12:
\end{verbatim}
We now try to solve this system of equations. After a few minutes of computing time
on a standard desktop PC,
Maple finds a symbolic representation (involving \verb|RootOf|s, which are  Maple 
representations for roots of polynomial equations)
of the general solution of the
system, for which we compute all possible values in numerical form. It
turns out that this computation (which again takes a few minutes) has to be done with very high
 precision, otherwise the results do not represent reasonable
 solutions with small residuals if substituted into the equations.
\begin{verbatim}
> sols12 := solve(eqs12, vars12):
> Digits := 200:
> FF := seq(evalf(allvalues(sol)), sol in sols12):
\end{verbatim}
We obtain 99 solutions altogether, 17 real solutions, and modulo complex conjugation $41$ different complex solutions.
Each  solution determines 8 parameters of the ansatz $S$. There remain  8 parameters to be determined compared
with 14 order conditions corresponding to the $14=22-8$ remaining
Lyndon words  over
$\nA=\{\AA_1,\AA_2,\AA_3,\AA_4\}$ of odd grade $\leq 8$.
It is remarkable that the resulting over-determined system of
 equations always has a solution.
To find a theoretical explanation for this fact
is the topic of current investigations. Here, however, it is verified by
a direct computation.

We select\footnote{
This is done for the purpose of presentation, the following considerations
apply to each of the 99 parameter sets. The solution selected for this
presentation leads to a particularly small local error.
%In the here displayed Maple code that solution involving only real numbers is
% selected which yields the smallest local error measure in the sense of
% \cite[Section~2.3]{hofiopmag}.
Note that the selected index 78 may belong to
different parameter sets in different runs of the code. }
one of the previously obtained 99 sets of 8 parameters, substitute it into
the ansatz $S$, and set up 4 equations corresponding to 4
(out of 9) selected Lyndon words  over $\nA$ of odd grade $\leq 8$ involving
$\AA_3$ but not $\AA_4$. The resulting system of equation is linear and
readily solved. That this solution solves also the equations corresponding to the $5=9-4$ not selected Lyndon words will be
verified below.
\begin{verbatim}
> F12 := FF[78]:
> W3 := [[A3], [A1, A1, A3], [A2, A3], [A1, A1, A1, A1, A3],
         [A1, A1, A2, A3], [A1, A1, A3, A2], [A1, A2, A1, A3],
         [A1, A3, A3], [A2, A2, A3]]:
> rhs3 := [0, 1/60, -1/30, 1/420, -1/168, 1/280, -1/840,
           1/420, -1/210]:
> vars3 := [f13, f23, f33, f43]:
> eqs3 := [seq(expand(wcoeff(w, subs(F12, S))), w in W3[2 .. 5])]
           - rhs3[2 .. 5]:
> F3 := op(solve(eqs3, vars3)):
\end{verbatim}
Analogously as before,
we substitute the 12 already obtained parameters into the ansatz
$S$, set up 4 equations corresponding to 4 (out of 5)
selected Lyndon words  over $\nA$ of odd grade $\leq 8$ involving
$\AA_4$, and solve the resulting linear system of equations. That the
obtained solution solves also the equation corresponding to the not
selected Lyndon word will again be verified below.
\begin{verbatim}
> W4 := [[A1, A4], [A1, A1, A1, A4], [A1, A2, A4],
         [A1, A4, A2], [A3, A4]]:
> rhs4 := [0, -1/840, 1/210, -1/140, -1/70]:
> vars4 := [f14, f24, f34, f44]:
> eqs4 := [seq(expand(wcoeff(w, subs(F12, F3, S))),
                w in W4[1 .. 4])] - rhs4[1 .. 4]:
> F4 := op(solve(eqs4, vars4)):
\end{verbatim}
Finally we print the calculated solution representing the 16
parameters $f_{j,k}$ of the ansatz $S$ and compute
its residual with respect to the order conditions corresponding to all
Lyndon words over $\nA$ of odd grade $\leq 8$ (including those not
previously selected). The tiny residual confirms that the
obtained scheme indeed satisfies the order conditions for order $p=8$ of
Theorem~\ref{Thm:order_conditions}.
\begin{verbatim}
> for y in [op(F12), op(F3), op(F4)] do
      lprint(evalf(y, 50))
  end do:
  f11 = -1.1210783473381738227756934594506597445892745485109
  f21 = 1.3210319274244662988569102191161576010502669814859
  f31 = -.11488794115695215928140654449977903918312514606917
  f41 = .41493436107065968320018978483428118272213271309425
  f12 = 1.0089705126043564404981241135055937701303470936598
  f22 = -1.1889339712738696420578749909323697235681087890036
  f32 = 0.44866039420480983666929215062389499923245100101695e-1
  f42 = -.13197275582656085011222031954705867101347489961070
  f13 = -.78475484313672167594298542161546182121249218395766
  f23 = .92477328275109744272940525314314765421496759253486
  f33 = 0.24950727790821017623386132247659342458740875944374e-1
  f43 = -.16496916740519678440980596377534517546121628452158
  f14 = .44843133893526952911027738378026389783570981940438
  f24 = -.52881775248948867348601923353730351864984279845615
  f34 = -0.24298790613584639672784191664606712944260031094723e-1
  f44 = .19795913373984127516833047932058800652021234941605
> W := [op(W12), op(W3), op(W4)]:
> RHS := [op(rhs12), op(rhs3), op(rhs4)]:
> printf("%.5e", max(map(abs, [seq(wcoeff(w,
         subs(F12, F3, F4, S)), w in W)]-RHS))):
  8.82689e-143
\end{verbatim}
%\begin{verbatim}
%> x := [1/2-sqrt((15+2*sqrt(30))*(1/140)),
%        1/2-sqrt((15-2*sqrt(30))*(1/140)),
%        1/2+sqrt((15-2*sqrt(30))*(1/140)),
%        1/2+sqrt((15+2*sqrt(30))*(1/140))]:
%> w :=  [1/4-(1/72)*sqrt(30), 1/4+(1/72)*sqrt(30),
%         1/4+(1/72)*sqrt(30), 1/4-(1/72)*sqrt(30)]:
%> P := [seq(unapply(add(
%        binomial(k, j)*binomial(k+j, j)*(-1)^(k+j)*X^j,
%        j = 0 .. k), X), k = 0 .. 4)]:
%> F := ArrayTools[Reshape](Matrix([[seq(op(2, f),
%       f in [op(F12),op(F3),op(F4)])]]), [4, 4]):
%> T := Matrix([seq([seq((2*l-1)*w[k]*P[l](x[k]), k = 1 .. 4)],
%       l = 1 .. 4)]):
%> evalf(F &* T), 20);
%\end{verbatim}
In the scheme
\begin{equation}\label{eq:cf_scheme}
S = \prod_{j=8,\ldots,1}\exp\big(\sum_{k=1}^4 f_{j,k}\AA_k\big)
\end{equation}
with parameters $f_{j,k}$ calculated by the above Maple code, the $\AA_k$
represent Legendre expansion coefficients defined by
integrals (\ref{eq:A_integral}). To obtain an effective numerical method we have
to approximate these integrals using a suitable quadrature formula.
Therefore we substitute
\begin{equation*}\label{eq:Al_substitution}
\AA_k\to(2k-1)\tau\sum_{l=1}^{K}w_kP_{k-1}(x_l)A(t_n+\tau x_l)
\end{equation*}
in (\ref{eq:cf_scheme})
with Gaussian nodes and weights
of order eight,
$$
(x_k) = \left(
\tfrac{1}{2}-\sqrt{\tfrac{15+2\sqrt{30}}{140}},\
\tfrac{1}{2}-\sqrt{\tfrac{15-2\sqrt{30}}{140}},\
\tfrac{1}{2}+\sqrt{\tfrac{15-2\sqrt{30}}{140}},\
\tfrac{1}{2}+\sqrt{\tfrac{15+2\sqrt{30}}{140}}
\right),
$$
$$
(w_k) =\left(
\tfrac{1}{4}-\tfrac{\sqrt{30}}{72},\
\tfrac{1}{4}+\tfrac{\sqrt{30}}{72},\
\tfrac{1}{4}+\tfrac{\sqrt{30}}{72},\
\tfrac{1}{4}-\tfrac{\sqrt{30}}{72}
\right),
$$
which corresponds to an application of Gaussian quadrature to
 (\ref{eq:A_integral}).
For the set of parameters $\{f_{j,k}\}$ displayed in the above
Maple code we obtain the integrator (cf.~(\ref{eq:non_auto_step}),
(\ref{eq:magnus_type_int}))
\begin{equation}\label{eq:cf_quad}
\nS(t_n,\tau) = \prod_{j=8,\ldots,1}\exp\big(\tau\sum_{k=1}^4 a_{j,k}A(t_n+\tau x_k)\big)
\end{equation}
with coefficients $a_{j,k}$ given in Table~\ref{tab:cf_coeffs_real}.
\begin{table}[t!]
\begin{center}
\caption{Coefficients $a_{j,k}$
for an 8th order commutator-free Magnus-type integrator (\ref{eq:cf_quad}).\label{tab:cf_coeffs_real}}
{\fontsize{6pt}{7pt}\selectfont
\begin{tabular}{cccc}
\hline
{\small $k=1$}&
{\small $k=2$}&
{\small $k=3$}&
{\small $k=4$}\\
\hline
\verb|-1.232611007291861933e+0|& \verb| 1.381999278877963415e-1|& \verb|-3.352921035850962622e-2|& \verb| 6.861942424401394962e-3| \\
\verb| 1.452637092757343214e+0|& \verb|-1.632549976033022450e-1|& \verb| 3.986114827352239259e-2|& \verb|-8.211316003097062961e-3| \\
\verb|-1.783965547974815151e-2|& \verb|-8.850494961553933912e-2|& \verb|-1.299159096777419811e-2|& \verb| 4.448254906109529464e-3| \\
\verb|-2.982838328015747208e-2|& \verb| 4.530735723950198008e-1|& \verb|-6.781322579940055086e-3|& \verb|-1.529505464262590422e-3| \\
\verb|-1.529505464262590422e-3|& \verb|-6.781322579940055086e-3|& \verb| 4.530735723950198008e-1|& \verb|-2.982838328015747208e-2| \\
\verb| 4.448254906109529464e-3|& \verb|-1.299159096777419811e-2|& \verb|-8.850494961553933912e-2|& \verb|-1.783965547974815151e-2| \\
\verb|-8.211316003097062961e-3|& \verb| 3.986114827352239259e-2|& \verb|-1.632549976033022450e-1|& \verb| 1.452637092757343214e+0| \\
\verb| 6.861942424401394962e-3|& \verb|-3.352921035850962622e-2|& \verb| 1.381999278877963415e-1|& \verb|-1.232611007291861933e+0| \\
 \hline
\end{tabular}
}
\end{center}
\end{table}
For these coefficients the positivity condition
\begin{equation}\label{eq:positivity_cond}
\mathrm{Re} f_{j,1}=\mathrm{Re}\sum_{k=1}^4a_{j,k} >0,\quad j=1,\dots,8,
\end{equation}
is not satisfied, in agreement with the fact that this cannot
 be the case for real coefficients,
see~\cite{positivitypaper2018}.
For some applications, however, it is essential for stability reasons
that this condition is satisfied.
This suggests to consider
schemes  with complex coefficients,
see~\cite{SergioFernandoMPaper2}. As was mentioned above,
of the 99 parameter sets $\{f_{j,k}\}$ which can be computed by
the above Maple code, $82=41\times 2$ involve complex numbers.
One of these parameter sets leads
to the scheme (\ref{eq:cf_quad}) with coefficients given
in Table~\ref{tab:cf_coeffs_complex},
for which (\ref{eq:positivity_cond}) is satisfied.

\begin{table}[t!]
\begin{center}
\caption{Real (top) and imaginary (bottom) parts of coefficients $a_{j,k}$
for an 8th order commutator-free Magnus-type integrator (\ref{eq:cf_quad}) satisfying the positivity condition (\ref{eq:positivity_cond}).\label{tab:cf_coeffs_complex}}
{\fontsize{6pt}{7pt}\selectfont
\begin{tabular}{cccc}
\hline
{\small $k=1$}&
{\small $k=2$}&
{\small $k=3$}&
{\small $k=4$}\\
 \hline
\verb| 5.162172083124911076e-2|& \verb|-5.787809823308952456e-3|& \verb| 1.404202563971892685e-3|& \verb|-2.873779919999358082e-4| \\
\verb| 1.129000600487386325e-1|& \verb|-1.811008163470541820e-2|& \verb| 8.982553129811831365e-3|& \verb|-2.544930699554437791e-3| \\
\verb| 2.631601314221973826e-2|& \verb| 1.983998701294184106e-1|& \verb|-4.965939955061425298e-2|& \verb| 1.197843408520720342e-2| \\
\verb|-1.592059248033346570e-2|& \verb| 1.424220211513735403e-1|& \verb| 4.842122146532602005e-2|& \verb|-1.013590436679991693e-2| \\
\verb|-1.013590436679991693e-2|& \verb| 4.842122146532602005e-2|& \verb| 1.424220211513735403e-1|& \verb|-1.592059248033346570e-2| \\
\verb| 1.197843408520720342e-2|& \verb|-4.965939955061425298e-2|& \verb| 1.983998701294184106e-1|& \verb| 2.631601314221973826e-2| \\
\verb|-2.544930699554437791e-3|& \verb| 8.982553129811831365e-3|& \verb|-1.811008163470541820e-2|& \verb| 1.129000600487386325e-1| \\
\verb|-2.873779919999358082e-4|& \verb| 1.404202563971892685e-3|& \verb|-5.787809823308952456e-3|& \verb| 5.162172083124911076e-2| \\
\hline
\verb|-1.187198036084005914e-1|& \verb| 1.331082409655082917e-2|& \verb|-3.229389682031679030e-3|& \verb| 6.609128526175740449e-4| \\
\verb| 1.359790143178213473e-1|& \verb| 3.226637801235380303e-3|& \verb|-5.647440118497178834e-3|& \verb| 1.831962429052182520e-3| \\
\verb|-1.952925932474600076e-2|& \verb| 4.339859420803126316e-2|& \verb| 4.884840043796339250e-3|& \verb|-1.849278537972746835e-3| \\
\verb| 3.513884130112852023e-3|& \verb|-7.185755041597012718e-2|& \verb| 1.591348406688517315e-2|& \verb|-1.887432258484616938e-3| \\
\verb|-1.887432258484616938e-3|& \verb| 1.591348406688517315e-2|& \verb|-7.185755041597012718e-2|& \verb| 3.513884130112852023e-3| \\
\verb|-1.849278537972746835e-3|& \verb| 4.884840043796339250e-3|& \verb| 4.339859420803126316e-2|& \verb|-1.952925932474600076e-2| \\
\verb| 1.831962429052182520e-3|& \verb|-5.647440118497178834e-3|& \verb| 3.226637801235380303e-3|& \verb| 1.359790143178213473e-1| \\
\verb| 6.609128526175740449e-4|& \verb|-3.229389682031679030e-3|& \verb| 1.331082409655082917e-2|& \verb|-1.187198036084005914e-1| \\
 \hline
\end{tabular}
}
\end{center}
\end{table}

\subsubsection*{Acknowledgements.} This work was supported in part by the Austrian
Science Fund (FWF) under grant P30819-N32 and the Vienna Science and Technology Fund (WWTF)
under grant MA14--002.

%\nocite{part1}

%\bibliographystyle{splncs04}
%\bibliography{casc2019}

\begin{thebibliography}{10}
\providecommand{\url}[1]{\texttt{#1}}
\providecommand{\urlprefix}{URL }
\providecommand{\doi}[1]{https://doi.org/#1}

\bibitem{alvfeh11}
Alverman, A., Fehske, H.: High-order commutator-free exponential
  time-propagation of driven quantum systems. J. Comput. Phys.  \textbf{230},
  5930--5956 (2011)

\bibitem{auzingeretal13c}
Auzinger, W., Herfort, W.: Local error structures and order conditions in terms
  of {L}ie elements for exponential splitting schemes. Opuscula Math.
  \textbf{34},  243--255 (2014)

\bibitem{casc2016}
Auzinger, W., Herfort, W., Hofst{\"a}tter, H., Koch, O.: Setup of order
  conditions for splitting methods. In: Gerdt, V.P., Koepf, W., Seiler, W.M.,
  Vorozhtsov, E.V. (eds.) Computer Algebra in Scientific Computing. pp. 30--42.
  Springer International Publishing, Cham (2016)

\bibitem{SergioFernandoMPaper2}
Blanes, S., Casas, F., Thalhammer, M.: High-order commutator-free
  quasi-{M}agnus exponential integrators for nonautonomous linear evolution
  equations. Comput. Phys. Commun.  \textbf{220},  243--262 (2017)

\bibitem{chin97}
Chin, S.: Symplectic integrators from composite operator factorizations.
  Physics Letters A  \textbf{226}(6),  344--348 (1997)

\bibitem{duval88}
Duval, J.: G{\'e}neration d'une section des classes de conjugaison et arbre des
  mots de {L}yndon de longueur born{\'e}e. Theoret. Comput. Sci.  \textbf{60},
  255--283 (1988)

\bibitem{expocon}
Expocon.mpl, \url{https://github.com/HaraldHofstaetter/Expocon.mpl}

\bibitem{haireretal02b}
Hairer, E., Lubich, C., Wanner, G.: Geometric Numerical Integration.
  Springer-Verlag, Berlin--Heidelberg--New York, 2nd edn. (2006)

\bibitem{hofiopmag}
{Hofst{\"a}tter}, H.: {Order conditions for exponential integrators} (Feb
  2019), submitted (for a preprint see arXiv:1902.11256).

\bibitem{positivitypaper2018}
Hofst{\"a}tter, H., Koch, O.: Non-satisfiability of a positivity condition for
  commutator-free exponential integrators of order higher than four. Numer.
  Math.  \textbf{141}(3),  681--691 (2019)

\bibitem{lothaire97}
Lothaire, M.: Combinatorics on Words. Encyclopedia of Mathematics and its
  Applications, Cambridge University Press, Cambridge, U.K. (1997)

\bibitem{MuntheKaas957}
Munthe{\textendash}Kaas, H., Owren, B.: Computations in a free {L}ie algebra.
  Phil. Trans. R. Soc. Lond. A  \textbf{357}(1754),  957--981 (1999)

\bibitem{suzuki95}
Suzuki, M.: New scheme of hybrid exponential product formulas with applications
  to quantum {M}onte--{C}arlo simulations. In: Computer Simulation Studies in
  Condensed-Matter Physics VIII. Springer-Verlag, Berlin, Heidelberg (1995)

\end{thebibliography}

\end{document}